\documentclass[12pt,oneside]{amsart}
\usepackage{geometry}                
\usepackage{graphicx}
\usepackage{amssymb}

\usepackage{amsmath,amsthm,amscd,amssymb}
\usepackage{latexsym}
\usepackage[colorlinks,citecolor=red,pagebackref,hypertexnames=false]{hyperref}
\numberwithin{equation}{section}
\theoremstyle{plain}
\newtheorem{theorem}{Theorem}[section]
\newtheorem{lemma}[theorem]{Lemma}
\newtheorem{corollary}[theorem]{Corollary}

\theoremstyle{definition}

\theoremstyle{remark}
\newtheorem{remark}[theorem]{Remark}

\newtheorem{case[theorem]}{Case}

\author{Alex Iosevich, Krystal Taylor and Ignacio Uriarte-Tuero}
\address{Department of Mathematics, University of Rochester, Rochester, NY}
\email{iosevich@math.rochester.edu}
\address{Department of Mathematics, The Ohio State, Columbus, OH}
\email{taylor.2952@osu.edu}
\address{Department of Mathematics, Michigan State University, MI}
\email{ignacio@math.msu.edu} 
\thanks{This work was partially supported by the NSA Grant H98230-15-1-0319 (A.I.) and grants DMS-1056965
(US NSF), MTM2010-16232, MTM2015-65792-P (MINECO, Spain) (I. U.-T.). }
\title{\parbox{14cm}{\centering{Pinned geometric configurations in Euclidean space and Riemannian manifolds}}}
\begin{document}
\maketitle
\begin{abstract} Let $M$ be a compact $d$-dimensional Riemannian manifold without a boundary. Given $E \subset M$, let $\Delta_{\rho}(E)=\{\rho(x,y): x,y \in E \}$, where $\rho$ is the Riemannian metric on $M$. Let $\Delta_{\rho}^x$ denote the pinned distance set, namely, $\{\rho(x,y): y \in E \}$ with $x \in E$. We prove that if the Hausdorff dimension of $E$ is greater than $\frac{d+1}{2}$, then there exist many $x \in E$ such that the Lebesgue measure of $\Delta^x_{\rho}(E)$ is positive. This result was previously established by Peres and Schlag in the Euclidean setting. The main result is deduced from a variable coefficient Euclidean formulation, which can be used to study a variety of geometric problems. We extend our result to the setting of chains studied in \cite{BIT15} and obtain a pinned estimate in this context. Moreover, we point out that our scheme is quite universal in nature and this idea will be exploited in variety of settings in the sequel. 
\end{abstract}  
\maketitle
\section{Introduction} 

The celebrated Falconer distance conjecture (see e.g. \cite{Fal86}, \cite{M95}) says that if the Hausdorff dimension of a compact set $E \subset {\Bbb R}^d$, $d \ge 2$, is greater than $\frac{d}{2}$, then the Lebesgue measure of the distance set $\Delta(E)=\{|x-y|: x,y \in E \}$ is positive. The best known results are due to Wolff in two dimension \cite{W04} and Erdogan \cite{Erd05} in higher dimensions. They proved that Lebesgue measure of the distance set is positive if the Hausdorff dimension of $E$ is greater than $\frac{d}{2}+\frac{1}{3}$. Under an additional assumption that $E \subset {\Bbb R}^2$ is Ahlfors-David regular, Orponen (\cite{O15}) proved that if the Hausdorff dimension of $E$ is greater than $1$, then the packing dimension of the distance set $\Delta(E)=\{|x-y|: x,y \in E \}$ is $1$, thus contributing significantly to the Falconer conjecture in this setting. 

An interesting variant of the Falconer distance problem is obtained by pinning the distance set. More precisely, given $x \in E$, let $\Delta^x(E)=\{|x-y|: y \in E\}$. Once again, the question is, how large does the Hausdorff dimension of $E$ need to be to ensure that the Lebesgue measure of $\Delta(E)$ is positive for $x$ in some subset of $E$. Peres and Schlag (\cite{PeSc2000}) proved that the conclusion holds if the Hausdorff dimension of $E \subset {\Bbb R}^d$, $d \ge 2$, is greater than $\frac{d+1}{2}$. Their method relies on a clever application of non-linear projection theory. 

In this paper we give a simple proof of the Peres-Schlag result and generalize it to a wide range of distance type functions. For $x \in E$, $E$ a compact subset of $\mathbb{R}^d$ for some $d\geq 2$, we consider
\begin{equation} \label{maindef} \Delta^x_{\phi}(E)=\{\phi(x,y): y \in E \},\end{equation} where $\phi: {\Bbb R}^d \times {\Bbb R}^d \to {\Bbb R}$ is continuous, infinitely differentiable, and satisfies 
\begin{equation}\label{nondeg}  
|\nabla_x\phi(x,y)|\neq 0, \ \  |\nabla_y\phi(x,y)|\neq 0 \end{equation}
almost everywhere on $E\times E$.  That is, for any Frostman measure $\mu$ on E, 
we set $F$ equal to the set of points where $\phi$ is not infinitely differentiable, and   
we assume that 
\begin{equation}\label{forbidden} \mu \times \mu\left( F \cup \left\{(x,y) :  |\nabla_y\phi(x,y)| =0 \text{  or  }  |\nabla_y\phi(x,y)|= 0) \right\}\right)=0.\end{equation}
\\
If, for instance, $\phi(x,y)=|x-y|$, then the gradient vanishes on the diagonal $x=y$.  
If $\phi(x,y)=|x-3y|$, then the gradient vanishes on $x=3y$.
If $\phi(x,y) = x\cdot y$, then the gradient is vanishes on the hyperplanes $x=0$ or $y=0$.  
\\

We also assume throughout that $\phi$ satisfies the non-vanishing Monge-Ampere determinant assumption: 
\begin{equation} \label{mongeampere} det
\begin{pmatrix} 
 0 & \nabla_{x}\phi \\
 -{(\nabla_{y}\phi)}^{T} & \frac{\partial^2 \phi}{dx_i dy_j}
\end{pmatrix} \end{equation} does not vanish on the set $\{(x,y): \phi(x,y)=t \}$, $t \not=0$. 

\vskip.125in 

Our main results are the following. 

\begin{theorem} \label{generaltech} Let $E$ be a compact subset of ${\Bbb R}^d$, $d \ge 2$ and let $\Delta^x_{\phi}(E)$ be defined as in (\ref{maindef}) above, with $\phi$ satisfying (\ref{nondeg}) and (\ref{mongeampere}). Suppose that the Hausdorff dimension of $E$ is greater than $\frac{d+1}{2}$. Then 
\begin{equation}\label{exceptionalset}dim_{{\mathcal H}} \left( \{x \in E: |\Delta^x_{\phi}(E)|=0 \} \right) \leq d+1-dim_{{\mathcal H}}(E).\end{equation}

In particular, if $dim_{{\mathcal H}}(E)=s_E>\frac{d+1}{2}$, then $|\Delta^x_{\phi}(E)|>0$ for $\mu$-a.e. $x \in E$, where $\mu$ is a Frostman measure on $E$. 
\end{theorem} 

\vskip.125in 

\begin{corollary} \label{riemanncor} Let $E$ be a closed subset of a compact $d$-dimensional, $d \ge 2$, Riemannian manifold $M$ without a boundary and let $\rho$ denote the Riemannian metric on $M$. Suppose that the Hausdorff dimension of $E$ is greater than $\frac{d+1}{2}$. Then  
\begin{equation} \label{riemannexceptionalset} dim_{{\mathcal H}} \left( \{x \in E: |\Delta^x_{\rho}(E)|=0 \} \right) \leq d+1-dim_{{\mathcal H}}(E). \end{equation}

In particular, if $dim_{{\mathcal H}}(E)=s_E>\frac{d+1}{2}$, then $|\Delta^x_{\rho}(E)|>0$ for $\mu$-a.e. $x \in E$, where $\mu$ is a Frostman measure on $E$. 
\end{corollary} 

Corollary \ref{riemanncor} follows from Theorem \ref{generaltech} and its proof by observing that in local coordinates, the Riemannian metric $\rho$ on a compact manifold without a boundary satisfies the assumptions (\ref{nondeg}) and (\ref{mongeampere}). See, for example, \cite{So93}, Chapter 6, and the references contained therein. 

\begin{remark} As the reader will see, the method of proof of Theorem \ref{generaltech} is very flexible. For example, assumption on the Monge-Ampere determinant can be easily replaced by a weaker condition on the rank of the Monge-Ampere matrix. More precisely, the only thing required to obtain a non-trivial result is $L^2({\Bbb R}^d) \to L^2_{\gamma}({\Bbb R}^d)$ Sobolev bound for some $\gamma>0$ for the generalized Radon transform 
$$ {\mathcal R}f(x)=\int_{\phi(x,y)=t} f(y) \psi(x,y) d\sigma_{x,t}(y),$$ where $\psi$ is a smooth cut-off function and $\sigma_{x,t}$ is a smooth surface measure on the set 
$$ \{y: \phi(x,y)=t\}.$$ 

It is also important to note that in this case of the Euclidean metric, $\phi(x,y)=|x-y|$, the only analytic input our proof uses is the fact that the operator 
$$ Af(x)=\int_{S^{d-1}} f(x-y) d\sigma(y),$$ where $\sigma$ is the Lebesgue measure on the sphere, maps $L^2({\Bbb R}^d))$ to 
$L^2_{\frac{d-1}{2}}({\Bbb R}^d)$. This is, of course, just equivalent to the well-known stationary phase estimate (see e.g. \cite{So93}) 
$$ |\widehat{\sigma}(\xi)| \leq C{|\xi|}^{-\frac{d-1}{2}}.$$ 
\end{remark} 

\vskip.125in 

\subsection{Pinned Chains}

Our methods also yield a result for pinned chains in thin subsets of Euclidean space.  
Let $E$ be a compact subset of $\mathbb{R}^d$, $d\geq 2$, and define 
a $k$-chain in $E$ with gaps $t= (t_1, \cdots, t_k)$ by 
$$\left\{ (x^1, \cdots, x^{k+1})\in E\times \cdots \times E:  \left(  \phi(x^1,x^2), \cdots, \phi(x^k,x^{k+1})  \right)= t \, \right\}.$$

Analogous to the distance set and pinned distance set, we define
\begin{equation}\label{chain}C_{k,\phi}(E)= \left\{   \left(  \phi(x,x^2), \cdots, \phi(x^k,x^{k+1})  \right): x, x^2, \dots, x^{k+1}\in E   \right\},\end{equation}
and
\begin{equation}\label{kchain}C_{k,\phi}^x(E)= \left\{   \left(  \phi(x,x^2), \cdots, \phi(x^k,x^{k+1})  \right): x^2, \dots, x^{k+1}\in E   \right\},\end{equation}
with $\phi$ satisfying (\ref{nondeg}) and (\ref{mongeampere}).  

\vskip.125in 

In the special case that $\phi$ is the Euclidean distance, it is shown in \cite{BIT15} that if the Hausdorff dimension of a set 
$E\subset \mathbb{R}^d$ is greater than $\frac{d+1}{2}$, then the $k$-dimensional Lebesgue measure
\begin{equation}\label{BIT} | C_{k,\phi}(E)|>0.\end{equation} 

Moreover, it is shown that under this assumption $C_{k,\phi}(E)$ has non-empty interior.  

\vskip.125in 

We now extend this result to show that $C_{k,\phi}^x(E)$ has positive $k$-dimensional Lebesgue measure for a large set of $x$ to be quantified.  

\vskip.125in

\begin{theorem}\label{pinnedtheoremchain}
Let $E$ be a compact subset of $\mathbb{R}^d$, $d\geq 2$ and let $C_{k,\phi}^x(E)$ as in \eqref{kchain}.
Suppose that the Hausdorff dimension of $E$ is greater than $\frac{d+1}{2}$.  
Then 
\begin{equation}\label{chainexceptionalset}\dim_{\mathcal{H}}  \left( \{  x\in E: \left|C_{k,\phi}^x(E)\right|=0  \} \right)  \le d+1 - \dim_{\mathcal{H}}(E).\end{equation}
In particular, if $\dim_{\mathcal{H}}(E)=s_E >\frac{d+1}{2}$, then $|C_{k,\phi}^x(E)| >0$ for $\mu$ a.e. $x\in E$, where $\mu$ is a Frostman measure on $E$.
\end{theorem}

\vskip.125in 

\begin{corollary}\label{pinnedtheoremmfldchain}
Let $E$ be a compact subset of a compact $d$-dimensional, $d \ge 2$, Riemannian manifold $M$ without a boundary and let $C_{k,\phi}^x(E)$ be defined as in \eqref{kchain}, where $\phi$ is replaced by $\rho$, the Riemannian metric on $M$. Suppose that the Hausdorff dimension of $E$ is greater than 
$\frac{d+1}{2}$. Then 
\begin{equation}\label{chainexceptionalsetmfld}\dim_{\mathcal{H}}  \left( \{  x\in E: \left|C_{k,\rho}^x(E)\right|=0  \} \right)  
\le d+1 - \dim_{\mathcal{H}}(E).\end{equation}
In particular, if $\dim_{\mathcal{H}}(E)=s_E >\frac{d+1}{2}$, then the $k$-dimensional Lebesgue measure of $C_{k,\rho}(E)$, denoted by 
$|C_{k,\rho}^x(E)| $ is positive for $\mu$ a.e. $x\in E$, where $\mu$ is a Frostman measure on $E$.
\end{corollary}

\vskip.125in 

Corollary \ref{pinnedtheoremmfldchain} follows from Theorem \ref{pinnedtheoremchain} in the same way Corollary \ref{riemanncor} follows from Theorem \ref{generaltech}. 

\vskip.125in 

\subsection{A general pinning scheme} While we do not study the pinning problem in full generality in this paper, we outline the basic general mechanism in this subsection and work out a few examples. For the sake of simplicity we work with the Euclidean distance, but the arguments work equally well with functions $\phi(x,y)$ satisfying (\ref{maindef}) and (\ref{nondeg}). Define the edge map 
$${\mathcal E}: \{1,2, \dots, k+1 \} \times \{1,2, \dots, k+1\} \to \{0,1\}$$ and let $n({\mathcal E})$ denote the number of non-zero values taken on by ${\mathcal E}$. 

Given a compact set $E \subset {\Bbb R}^d$, a positive integer $k \ge 1$ and an edge map ${\mathcal E}$, we define the $k$-point configuration 
${\mathcal P}_{k, {\mathcal E}}(E)$ to be the set of $n=n({\mathcal E})$ vectors with entries $|x^i-x^j|$ where $1 \leq i<j \leq k+1$, 
${\mathcal E}(i,j)=1$. In this way, we may naturally view ${\mathcal P}_{k,{\mathcal E}}(E) \subset {\Bbb R}^{n({\mathcal E})}$. For example, the chain set from the previous subsection corresponds to the edge function ${\mathcal E}(i,j)=1$, $j=i+1$, 
$1 \leq i \leq k$, and $0$ otherwise. 

To illustrate the general pinning process, we begin with a $2$-chain and pin the middle vertex. More precisely, we have 
$$ {\mathcal E}(1,2)={\mathcal E}(2,3)=1$$ and the rest of the values of ${\mathcal E}(i,j)=0$. Since we are pinning the middle vertex, we are looking at the set 
$$ {\mathcal P}^x_{2, {\mathcal E}}(E)=\{(|x^1-x|, |x^2-x|): x^1, x^2 \in E\}.$$ 

Following the argument for pinned chains below in this context, and relabelling the coordinates, we see that proving that the two-dimensional Lebesgue measure of 
${\mathcal P}^x_{2, {\mathcal E}}(E)$ is positive, we must estimate 
\begin{equation} \label{L4} \int \epsilon^{-4} \mu \times \mu \times \mu \times \mu \times \mu 
\{(x^1,x^2,x^3,x^4,x^5) \in E^5: t_i-\epsilon \leq |x^i-x^3| \leq t_i+\epsilon; i \not=3 \} dt. \end{equation}  

To put it another way, after relabelling the coordinates, we are led to consider ${\mathcal P}_{4, {\mathcal E}'}(E)$, where ${\mathcal E}'(i,j)=1$ if $1 \leq i \leq 4$ and $j=5$, and $0$ otherwise. This is a four vertex star with the center at $x^5$. This instantly leads to an interesting estimate because estimation of (\ref{L4}) leads to an $L^4$ estimate of the underlying operator (see (\ref{opbound}) below), instead of the $L^2$ bound required to prove Theorem \ref{generaltech} and Theorem \ref{pinnedtheoremchain}. In this case integration in $t$ is helpful, in view of Mockenhaupt-Seeger-Sogge local smoothing estimates. 

In general, in order to study the pinned version of ${\mathcal P}_{k, {\mathcal E}}(E)$, where the $k+1$st vertex is pinned (which we can always arrange by relabeling), we are led to consider ${\mathcal P}_{2k, {\mathcal E}'}(E)$, where ${\mathcal E}'(i,j)={\mathcal E}(i,j)$ if $1 \leq i<j \leq k+1$, ${\mathcal E}'(k+1, k+1+j)={\mathcal E}(j,k+1)$, ${\mathcal E}'(k+1+i, k+1+j)={\mathcal E}(i,j)$ and ${\mathcal E}'(i,j)=0$ if $i<k+1$ and $j>k+1$. 

\begin{figure}
\label{twonecklaces-figure}
\centering
\includegraphics[scale=.5]{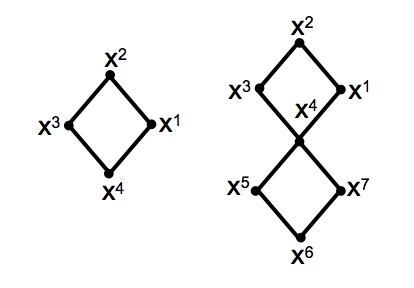}
\caption{Pinned necklace reduces to the consideration of two necklaces sharing the pinned point}
\end{figure}

\vskip.125in 

A slightly more complicated situation is described in Figure 1 above. In this case ${\mathcal E}(1,2)={\mathcal E}(2,3)={\mathcal E}(3,4)= {\mathcal E}(1,4) = 1$ and ${\mathcal E}(i,j)=0$ otherwise. Denote this configuration by ${\mathcal P}_{4,{\mathcal E}}$, which is precisely what is depicted on the left hand side of Figure 1. If we pin the vertex $x^4$ and apply our method, we arrive at the configuration ${\mathcal P}_{2k, {\mathcal E}'}$ depicted on the right hand side of Figure 1. 

Our method also allows us to pin two or more vertices. We shall undertake a systematic study of the pinned configurations in a subsequent paper. 

\vskip.25in 

\section{Proof of Theorem \ref{generaltech}} 

\subsection{Basic reductions}\label{basicreductions} Let $\mu$ be a Frostman measure supported on $E$. Recall that this means that 
\begin{equation} \label{muball} \mu(B(x,r)) \leq Cr^{s_{\mu}} \end{equation} for $s_{\mu}<dim_{{\mathcal H}}(E)$ and we can make $s_{\mu}$ arbitrarily close to $dim_{{\mathcal H}}(E)$. See, for example, \cite{M95} for the construction. 

Define the measure $\nu_x$ on $\Delta^x_{\phi}(E)$ by the relation 
\begin{equation} \label{measdef} \int f(t) d\nu_x(t)=\int f(\phi(x,y)) d\mu(y). \end{equation} 

Applying Cauchy-Schwartz we see that (formally) 
\begin{equation} \label{csformal} 1={\left( \int d\nu_x(t) \right)}^2 \leq |\Delta^x_{\phi}(E)| \cdot \int \nu_x^2(t) dt. \end{equation} 

If we could make sense of the right hand side and show that 
$$ \int \int \nu_x^2(t) d\mu(x) dt \leq C,$$ we could conclude that 
$$ \int \nu_x^2(t) dt<\infty$$ for $\mu$-a.e. $x \in E$, and plugging this into (\ref{csformal}) would show that $|\Delta^x_{\phi}(E)|>0$ for $\mu$-a.e. $x \in E$. 

In order to prove the more precise estimate (\ref{exceptionalset}) we will show that if $\lambda$ is a compactly supported Borel measure such that 
\begin{equation} \label{lambdaball} \lambda(B(x,r)) \leq Cr^{s_{\lambda}} \end{equation} for some $s_{\lambda}>0$, then 
\begin{equation} \label{L2BoundOnNuXInDLambdaDT} \int \int \nu_x^2(t) d\lambda(x) dt \leq C  \end{equation} if 
$$s_{\lambda}+dim_{{\mathcal H}}(E)>d+1.$$ 

The conclusion (\ref{exceptionalset}) is recovered by taking $\lambda$ to be a Frostman measure on a subset of $E$. 

\vskip.125in 

To make the setup rigorous and in order to understand the geometric meaning of the quantity 
$$ \int \int \nu_x^2(t) d\lambda(x) dt,$$ let $\rho$ be a smooth cut-off function, supported in the ball of radius $2$ and identically equal to $1$ in the ball of radius $1$ centered at the origin, with 
$\int \rho=1$. (We used the same notation $\rho$ for the Riemannian metric, but the context makes it clear which is which.) Let $\rho_{\epsilon}(u)=\epsilon^{-1} \rho(u/\epsilon)$. Then by (\ref{measdef}), 
$$ \nu_x*\rho_{\epsilon}(t)=\epsilon^{-1} \int \rho \left(\frac{t-u}{\epsilon} \right) d\nu_x(u)$$ 
\begin{equation}\label{keyq} = \epsilon^{-1} \int  \rho \left(\frac{t-\phi(x,y)}{\epsilon} \right) d\mu(y).\end{equation}

Note that the right hand side is clearly in $L^2(\lambda)$ with constants possibly depending on $\epsilon$. Our goal is to prove that these constants are, in fact, independent of $\epsilon$. 

Also note that 
$$ \int \nu_x*\rho_{\epsilon}(t) dt=\epsilon^{-1} \int \int \rho \left( \frac{t-u}{\epsilon} \right) d\nu_x(u)dt$$ 
\begin{equation} \label{stillprob}=\int d\nu_x(u)=1. \end{equation} 

Since the integration in $t$ is compact owing to the compactness of $E$, we can insert a smooth cut-off function $\beta$.  Squaring  and integrating with respect to $d\lambda(x)$, we obtain 
\begin{equation}\label{LHS} \int \int {\left| \frac{1}{\epsilon}  \int  \rho \left(\frac{t-\phi(x,y)}{\epsilon} \right) d\mu(y) \right|}^2 d\lambda(x) \beta(t)dt\end{equation}
$$=\int \int \frac{1}{\epsilon^2} \int \int \rho \left(\frac{t-\phi(x,y)}{\epsilon} \right) \rho \left(\frac{t-\phi(x,z)}{\epsilon} \right) d\mu(y) 
d\mu(z) d\lambda(x) \beta(t)dt.$$ 

\vskip.125in 

In view of (\ref{stillprob}) and the formal argument in (\ref{csformal}), with $\nu_x$ replaced by $\nu_x^{\epsilon} = \nu_x*\rho_{\epsilon}$ (see \cite{Mat2015}, pg.17), it suffices for us to show that 

\begin{equation} \label{hingeappears} \frac{1}{\epsilon^2} \int \lambda \times \mu \times \mu \{(x,y,z): t-\epsilon \leq \phi(x,y) \leq t+\epsilon; \ t-\epsilon \leq \phi(x,z) \leq t+\epsilon \} \beta(t) dt \leq C \end{equation}  independently of $\epsilon$. This will give us a lower bound on the Lebesgue measure of the support of $\nu_x^{\epsilon}$ which is independent of $\epsilon$. As an immediate consequence, we deduce the lower bound on the Lebesgue measure of $\Delta_{\phi}^x(E)$. 
\\

To be more precise, we consider \eqref{measdef}, but with the measures convolved with $\rho_{\epsilon}$. Now $\nu_x^{\epsilon}$ is indeed absolutely continuous with respect to Lebesgue measure. Then we prove \eqref{L2BoundOnNuXInDLambdaDT} for $\nu_x^{\epsilon}$. Note that $\nu_x^{\epsilon}$ converges to $\nu_x$ as measures and hence as distributions. Moreover, $\parallel  \nu_x^{\epsilon}  \parallel_{L^2(d\lambda \times dt)} \leq C$, so by Banach-Alaouglu we get that $\nu_x$ is absolutely continuous with respect to Lebesgue measure and has density in $L^2(d\lambda \times dt)$. Now \eqref{csformal} is indeed rigorous.
\\

The expression on the left hand side of (\ref{hingeappears}), excluding the integration in $t$, is precisely the quantity that arises in the study of geometric hinges, first explored by the second listed author in her thesis and later studied systematically in the case $\phi(x,y)=|x-y|$ in \cite{BIT15} and \cite{IKSTU16}. 
\\

We have shown that the proof of Theorem \ref{generaltech} would follow if we could establish the estimate (\ref{hingeappears}). 
We now reduce this estimate to an operator bound. Let $\psi$ be a smooth cut-off function which is equal to $1$ on $\left(E\times E\right)\backslash N$, where 
$N$ is a small open neighborhood of the set of points being measured in equation \eqref{forbidden}. By the regularity of the measure $\mu\times \mu$ (see see \cite{Ru87},theorem 2.18) and by assumption \eqref{forbidden}, we can choose $N$ so that $\mu\times \mu(N)$ is arbitrarily small.  Further, assume that $\psi$ vanishes on $N$.   Set

\begin{equation}\label{operatordefn} T_{\phi, t}^{\epsilon}f(x)=\frac{1}{ \epsilon} \int  \rho\left(\frac{t-\phi(x,y)}{\epsilon}\right) f(y) \psi(x,y) dy,\end{equation} where $\rho$ was introduced above equation \eqref{keyq} and $\phi$ is as in the statement of the theorem.
\\

The estimate (\ref{hingeappears}) would follow instantly from the operator bound
\begin{equation} \label{opbound} {\left( \int {|T_{\phi}^{\epsilon}(f\mu)(x)|}^2 d\lambda(x) \right)}^{\frac{1}{2}} \leq C {||f||}_{L^2(\mu)}
\end{equation} under the assumption $s_{\lambda}>d+1-s_{\mu}$. Here, the constant $C$ is uniform in $t$ and independent of $\epsilon$, and $T_{\phi}^{\epsilon}(f\mu)(x)$ is the supremum in $t > 0$ of $ T_{\phi, t}^{\epsilon}(f\mu)(x)$ .  This is where we now turn our attention. 

\vskip.125in 

\subsection{Proof of the operator bound (estimate (\ref{opbound}))} 

We prove estimate \eqref{opbound} by showing that, for $g\in L^2(\lambda)$, 
$$ \lvert <T_{\phi}^{\epsilon}f\mu, g\lambda> \rvert \leq C{||f||}_{L^2(\mu)} {||g||}_{L^2(\lambda)}.$$ 
\vskip.125in 

Let $\alpha_0(\xi)$ and $\alpha$ be smooth cut-off functions such that $\alpha_0$ is supported in the ball $\{|\xi|<4\}$, $\alpha$ is supported in the annulus $\{1/2\le |\xi| \le 4\}$, and $\alpha_0(\xi) + \sum_k \alpha(2^{-j} \cdot) \equiv 1$.

Define $P_jT_{\phi}^{\epsilon}$, the classical Littlewood-Paley projection, by the relation 
$$ \widehat{P_jT_{\phi}^{\epsilon}}=\widehat{T_{\phi}^{\epsilon}} \cdot \alpha(2^{-j} \cdot),$$
and in general, for $f \in L^2(dx)$,
$$ \widehat{P_jf}=\widehat{f} \cdot \alpha(2^{-j} \cdot).$$

Let $g \in L^2(\lambda)$. Then 
\begin{equation}\label{TEpsilonPhiIscontinuousInL2AndHenceCanTakeInfinteSumOfLPProjectionsOut} \left<T_{\phi}^{\epsilon}f\mu, g\lambda \right>=\sum_{j,k\geq  0} \left<T_{\phi}^{\epsilon}(P_jf\mu), P_kg\lambda)\right>.\end{equation}

This equation would probably benefit from a few words of explanation. Indeed, $T_{\phi}^{\epsilon}$ is trivially continuous $L^2(dx) \to L^2(dx)$, and thus, in the $L^2(dx)$ sense, $T_{\phi}^{\epsilon}f\mu = \sum_{j\geq  0} T_{\phi}^{\epsilon}(P_jf\mu)$, if $f\mu$ were in $ L^2(dx)$, and then $ \left<\cdot , \cdot \right>$ denotes inner product in $L^2(dx)$, so the $\sum_{j,k\geq  0}$ can be taken out of the inner product. This is indeed what happens, except that one has to make sense of $f\mu$ and $g\lambda$, which are not in $ L^2(dx)$, since the measures $\mu$ and $\lambda$ are typically singular with respect to Lebesgue measure.  

In order to do that, instead of $\mu$ and $\lambda$, consider $\mu_\theta = \mu * \rho_{\theta}$ and $\lambda_\delta = \lambda * \rho_{\delta}$ (since $\epsilon$ has already been used as a parameter for regularization). Then these new measures are absolutely continuous with respect to Lebesgue measure and we can consider $f \in L^2(\mu_\theta)$ and $g \in L^2(\lambda_\delta)$ which are continuous and compactly supported. Now we can run the whole argument of the proof below, and get \eqref{opbound} for the measures $\mu_\theta$ and $\lambda_\delta$. Now note that by Fubini, for any $F$ continuous and compactly supported,

\begin{equation}\label{PassingMollificationFromMeasureToFunction}
\int F(x) \; d\mu_\theta(x) = \int F(x) \left[ \int \frac{1}{\theta} \ \rho \left( \frac{x - u}{\theta} \right) d\mu (u) \right] \; dx = \int \left( F * \rho_\theta  \right) (u) \; d\mu (u).
\end{equation}

Recall that if $F$ is continuous and compactly supported, $F * \rho_\theta \to F$ uniformly. Thus, by the dominated convergence theorem we have that, as $\theta \to 0$,

\begin{equation}\label{ConvergenceOfL2NormsOfMollifiedMeasuresToL2NormOfUnmollifiedMeasure}
\parallel  f \parallel^2_{L^2(\mu_\theta)} = \int |f|^2 * \rho_\theta \ d\mu   \to  \int |f|^2 \ d\mu = \parallel f\parallel^2_{L^2(\mu)}. 
\end{equation}

On the other hand we claim that, for all $x$, as $\theta \to 0$, $ T_{\phi}^{\epsilon}(f\mu_\theta)(x)    \to T_{\phi}^{\epsilon}(f\mu)(x)$. Indeed, fix $x$ and note that $\psi_{x}(y) = \psi(x,y)$ is then a continuous function in $y$ with compact support. Then the same reasoning as \eqref{PassingMollificationFromMeasureToFunction} and \eqref{ConvergenceOfL2NormsOfMollifiedMeasuresToL2NormOfUnmollifiedMeasure} shows that as $\theta \to 0$, 
\begin{eqnarray}
T_{\phi}^{\epsilon}(f\mu_\theta)(x) =   \frac{1}{ \epsilon} \int  \rho\left(\frac{t-\phi(x,y)}{\epsilon}\right) f(y) \psi(x,y) d\mu_\theta(y) = \nonumber \\
 \frac{1}{ \epsilon} \int  \left[ \rho\left(\frac{t-\phi(x,\cdot)}{\epsilon}\right) f(\cdot) \psi_{x}(\cdot) \right]_{\theta}(u)  d\mu(u) \to 
T_{\phi}^{\epsilon}(f\mu)(x). \nonumber
\end{eqnarray}

As a consequence, by Fatou's lemma, the same reasoning as \eqref{PassingMollificationFromMeasureToFunction} and \eqref{ConvergenceOfL2NormsOfMollifiedMeasuresToL2NormOfUnmollifiedMeasure} for $\lambda_\delta$, and duality, we get 
\eqref{opbound} (for the measures $\mu$ and $\lambda$, as stated), and for $f$ continuous and compactly supported, which, by density, yields \eqref{opbound} in full generality.

Returning to \eqref{TEpsilonPhiIscontinuousInL2AndHenceCanTakeInfinteSumOfLPProjectionsOut}, we consider the sum over $|j-k|> K$, for a large integer $K$ to be determined, and over $|j-k|\le K$ separately.

\begin{lemma} \label{local} Let $K$ be a positive integer.  Then 
\begin{equation}\label{localeq} \sum_{|j-k| \leq K}\int \widehat{T_{\phi}^{\epsilon}(P_jf\mu)}(\xi) 
\widehat{P_kg\lambda}(\xi) d\xi \le 
C{||f||}_{L^2(\mu)} {||g||}_{L^2(\lambda)}.
\end{equation}
\end{lemma} 

\begin{lemma} \label{microlocal} With the notation above, there exists a positive integer $K$ such that
\begin{equation}\label{microlocaleq} \sum_{|j-k| >K}\int \widehat{T_{\phi}^{\epsilon}(P_jf\mu)}(\xi) 
\widehat{P_kg\lambda}(\xi) d\xi \le 
C{||f||}_{L^2(\mu)} {||g||}_{L^2(\lambda)}.
\end{equation}
\end{lemma} 
\vskip.125in 

To prove Lemma \ref{local}, we note that since the support of $P_k$ in the frequency side is contained in that of $P_{k-1} + P_k + P_{k+1}$, we can add these in front of the left hand term:
\begin{gather}
 \sum_{|j-k|\le K} \left< T_{\phi}^{\epsilon}(P_jf\mu) , P_kg\lambda)\right> = \sum_{|j-k|\le K} \left< P_{k-1} T_{\phi}^{\epsilon}(P_jf\mu)      , P_kg\lambda)\right> +  \nonumber \\
\sum_{|j-k|\le K} \left< P_k T_{\phi}^{\epsilon}(P_jf\mu)      , P_kg\lambda)\right> + \sum_{|j-k|\le K} \left< P_{k+1} T_{\phi}^{\epsilon}(P_jf\mu)      , P_kg\lambda)\right> \nonumber
\end{gather}

We only consider the second term above, the first and the third being similar. Then we have that

\begin{align} \label{almost} 
\sum_{|j-k|\le K} & \left< P_k T_{\phi}^{\epsilon}(P_jf\mu)      , P_kg\lambda)\right> \leq \nonumber \\
& \sum_{|j-k| \leq K}   
{\left(\int {|\widehat{P_kg\lambda}(\xi)|}^2 {|\xi|}^{-d+s_{\lambda}-\eta} d\xi \right)}^{-\frac{1}{2}}  
\cdot
{\left(\int |    ( P_k T_{\phi}^{\epsilon}(P_j f\mu))^{\widehat{}}  \,    |^2 {|\xi|}^{d-s_{\lambda}+\eta} d\xi \right)}^{\frac{1}{2}},
\end{align}
for $\eta > 0$ sufficiently small.

We shall need the following lemma from \cite{BIT15} (Lemma 2.5). We shall give a proof below for the sake of completeness. 
\begin{lemma} \label{fenergy} Let $\lambda$ be a compactly supported Borel measure such that $\lambda(B(x,r)) \leq Cr^{s_{\lambda}}$ for some $s_{\lambda} \in (0,d)$. Suppose that $\gamma>d-s_{\lambda}$. Then for $g \in L^2(\lambda)$, 
\begin{equation} \label{fenergyeq} \int {|\widehat{g\lambda}(\xi)|}^2 {|\xi|}^{-\gamma} d\xi \leq C'{||g||}^2_{L^2(\lambda)}. \end{equation}
\end{lemma} 

Applying Lemma (\ref{fenergy}) shows that the first term on the right hand side in \eqref{almost} is uniformly bounded. To handle the second term, we shall need the following result due to Phong and Stein (\cite{PhSt91}). 
\begin{theorem} \label{phongsteinth} Let $T_{\phi}^{\epsilon}$ be defined as above with $\phi$ satisfying assumptions (\ref{nondeg}) and (\ref{mongeampere}). Then 
$$ T_{\phi}^{\epsilon}: L^2({\Bbb R}^d) \to L^2_{\frac{d-1}{2}}({\Bbb R}^d) \ \text{with constants independent of} \ \epsilon,$$ where $L^2_{\gamma}({\Bbb R}^d)$ denotes the Sobolev space of function with $\gamma$ (generalized) derivatives in $L^2({\Bbb R}^d)$. 
\end{theorem} 

The presence of $P_k$ implies there exist constants $0<c\le C $ such that the square of the second term in the summand of (\ref{almost}) is bounded by 
$$ \int {|  ( P_k T_{\phi}^{\epsilon}(P_jf\mu)  )^{\widehat{}}(\xi)|}^2 {|\xi|}^{d-s_{\lambda}+\eta} d\xi 
\leq C 2^{k(d-s_{\lambda}+\eta)} \int_{c2^j \leq |\xi| \leq C2^j} {|(   T_{\phi}^{\epsilon}(P_j f\mu)   )^{\widehat{}}    (\xi)|}^2 d\xi.$$

By Theorem \ref{phongsteinth}, and recalling that $|j-k|\le K$, the right hand side is bounded by 
\begin{equation} \label{arithmetic} C2^{j(d-s_{\lambda}+\eta)} \cdot 2^{-j(d-1)} \int {|\widehat{P_jf\mu}(\xi)|}^2 d\xi. \end{equation} 

It easily follows from Lemma (\ref{fenergy}) that 
$$  \int {\left|\widehat{P_jf\mu}(\xi)\right|}^2 d\xi \leq C 2^{j(d-s_{\mu}+\eta')} {||f||}^2_{L^2(\mu)},$$
for $\eta' > 0$ sufficiently small.

Inserting this back into (\ref{arithmetic}) we obtain 
$$ C2^{j(d-s_{\lambda}+\eta)} \cdot 2^{-j(d-1)} \cdot 2^{j(d-s_{\mu}+\eta')} {||f||}^2_{L^2(\mu)}=C2^{j(d+1-s_{\mu}-s_{\lambda}+\eta+\eta')} {||f||}^2_{L^2(\mu)}.$$

The geometric series converges precisely when $s_{\lambda}>d+1-s_{\mu}$ thus completing the proof of Lemma \ref{local}.  
\\

To complete the proof of Theorem \ref{generaltech}, it remains to prove Lemma \ref{microlocal} and Lemma \ref{fenergy}. These proofs can be found in the next two subsections. 

\vskip.125in 

\subsection{Proof of Lemma \ref{microlocal}} (See, for example, \cite{EIT11}, for similar arguments). 
Use Fourier inversion to write 
\begin{align*}
T_{\phi}^{\epsilon}(P_j f\mu) (x)
&= \int \Phi_{\epsilon}(\phi(x,y) -t)   P_j(f\mu)(y) \psi(x,y) dy\\
&= \iiint \widehat{\Phi_{\epsilon}}(s)     \widehat{ P_j(f\mu)}(\zeta)  \,   \exp\left(   2\pi\, i(   (\phi(x,y) -t)s  + y\cdot\zeta )    \right) \psi(x,y) ds d\zeta  dy\\ 
\end{align*}
where $\Phi_{\epsilon} $ is $\frac{1}{2\epsilon}$ times the indicator function of the interval $[-\epsilon, \epsilon]$.  \\

Recalling that $\left(P_j f\mu)\right)^{\widehat{}}  (\zeta) = \widehat{ f\mu}(\zeta)   \alpha(2^{-j}\zeta) $, we have
\begin{equation}\label{name}
 \left(T_{\phi}^{\epsilon}(P_j f\mu)\right)^{\widehat{}}  (\xi)
= \iint \widehat{\Phi_{\epsilon}}(s)   \,  \widehat{ f\mu}(\zeta)  \,  \alpha(2^{-j}\zeta) \,G(s,\xi, \zeta) ds d\zeta,
\end{equation}
where $G(s, \xi,\zeta) = \iint  e^{\left(   2\pi\, i(   (\phi(x,y) -t)s  + y\cdot\zeta  - x\cdot \xi  )    \right)} \psi(x,y) dx dy.$\\
\\

We now give an upper bound on the modulus of $G$ in the regime that $|\zeta|\sim 2^j$ and $|\xi| \sim 2^k$.
\begin{lemma}\label{IBP} Suppose that $|\zeta|\sim 2^j$ and $|\xi| \sim 2^k$.  Then there exists a $K>0$ so that if $|j-k|>K$, then for each positive integer $M$, there exist a positive constants $c_M>0$ so that 
$$|G(s,\xi,\zeta)|  \le c_M\, \inf \left\{|s|^{-M}, 2^{-j\,M}, 2^{-k\,M} \right\} \le c_{\tilde{M}}\, \inf \left\{ |s|^{-M}, 2^{-j\,(d-s_{\lambda})/2}, 2^{-k\,(d-s_{\mu})/2}   \right\} . $$ \end{lemma}

We give the proof of this Lemma momentarily. Returning to the calculation above, we multiply both sides of \eqref{name} by $ \alpha(2^{-k}\xi)$ to see that 
\begin{equation}\label{boundme}  
 \alpha(2^{-k}\xi)  \left( T_{\phi}^{\epsilon}(P_j f\mu) \right)^{\widehat{}}(\xi)
=   \iint  \left( \Phi_{\epsilon}  \right)^{\widehat{}}(s)     \widehat{ (f\mu)}(\zeta)  \,  \alpha(2^{-k}\xi) \alpha(2^{-j}\zeta)G(s,\xi,\zeta) ds d\zeta.\end{equation}

Inserting the estimate from Lemma \ref{IBP} and integrating in $s$,
we bound \eqref{boundme} above in absolute value by 
$$ \alpha(2^{-k}\xi)  \min\{ 2^{-j(M-1)} ,2^{-k(M-1)} \}  \int  \left| \widehat{(f\mu)}(\zeta)\right|  \,  \alpha(2^{-j}\zeta) d\zeta.$$

After applying Cauchy Schwarz, we use Lemma \ref{fenergy} to bound this expression above by
$$ \alpha(2^{-k}\xi)  \min\{ 2^{-j(M-1)} ,2^{-k(M-1)} \} 2^{j(2d-s_{\mu}+\eta'')/2} \|f\|_{L^2(\mu)} ,$$
for $\eta'' > 0$ sufficiently small.

At last, 
\begin{align*}
\lvert <  & \left(P_k T_{\phi}^{\epsilon}(P_j f\mu)  \right)^{\widehat{}},  \widehat{g\lambda}> \rvert \lesssim  \\
& \int \alpha(2^{-k}\xi)  \left   |\widehat{g\lambda}(\xi)   \right| d\xi  \min\{ 2^{-j(M-1)} ,2^{-k(M-1)} \} 2^{j(2d-s_{\mu}+\eta'')/2} \|f\|_{L^2(\mu)}  \\
&\lesssim    \min\{ 2^{-j(M-1)} ,2^{-k(M-1)} \} 2^{j(2d-s_{\mu}+\eta'')/2} \|f\|_{L^2(\mu)}  2^{k(2d-s_{\lambda}+\eta''')} \| g\|_{L^2(\lambda)},\end{align*}
for $\eta''' > 0$ sufficiently small.

Summing ${\displaystyle{\sum_{j,k: |j-k| >K} }}$ in the previous display, yields convergent geometric series, thus proving the lemma. 


\vskip.125in 

\subsection{Proof of Lemma \ref{IBP}} We prove the first inequality in the statement of the Lemma, and the second follows by taking $M$ sufficiently large. 
We first compute the critical points of the phase function for 
 \begin{equation}\label{G}G(s, \xi,\zeta) = \iint  e^{\left(   2\pi\, i(   (\phi(x,y) -t)s  + y\cdot\zeta  - x\cdot \xi  )    \right)} \psi(x,y) dx dy.\end{equation}

The critical points occur when 
\begin{equation}\label{critical} s\nabla_x\phi(x,y) -\xi =0 \hspace{.2in}  \text{ and } \hspace{.2in}   s\nabla_y\phi(x,y) +\zeta= 0.
\end{equation}
We may assume that exists positive constants $0<c<C$ so that 
$$c<\left| \nabla_x\phi(x,y) \right| <C    \hspace{.2in}  \text{ and } \hspace{.2in}   c<\left| \nabla_y\phi(x,y) \right| <C  ,$$
where the lower bound follows from the non degeneracy assumption in \eqref{nondeg} coupled with the assumption that $\psi$ has compact support. \\

We argue that if $|\zeta|\sim 2^j$, $|\xi| \sim 2^k$, and $j$ and $k$ are sufficiently separated, then both equations in \eqref{critical} cannot both simultaneously hold for $(x,y)$ in the support of $\psi$.  Let us first take care of the scenario when $s=0$.  
If $s=0$, then we may at least assume that either $|\xi|\neq 0$ or $|\zeta| \neq 0$, and the lemma follows by repeated integration by parts.
Assume then that $s\neq0$.\\

If both equations in \eqref{critical} were to hold, then it would follow that $$\frac{|  \xi |}{|  \zeta |} \in \left( \frac{c}{C}, \frac{C}{c}  \right).$$
This is clearly false if $|\zeta|\sim 2^j$, $|\xi| \sim 2^k$, and $j$ and $k$ are sufficiently separated. We conclude then that there exists $K>0$ so that if $|j-k|>K$, then the integrand of $G$ is supported away from critical points.\\

Assuming then that $|j-k|>K$, it follows that for each $(x,y)$ in the support of $\psi$ that either 
\begin{equation}\label{away} s\nabla_x\phi(x,y) -\xi \neq 0 \hspace{.2in}  \text{ or } \hspace{.2in}   s\nabla_y\phi(x,y) +\zeta \neq 0.
\end{equation}

Fix $(x,y)$ in the support of $\psi$, and assume the first inequality in \eqref{away} holds.  
Let $l\ \in \left\{1, \cdots, d\right\}$ denote the coordinate such that $| s\frac{ \partial \phi(x,y)}{\partial x_l} -\xi_l |$  is maximal.  
It follows that $| s\frac{ \partial \phi(x,y)}{\partial x_l} -\xi_l |   \geq    \frac{1}{d} |   s\nabla_x\phi(x,y) -\xi  |$.  
We see that $ e^{\left(   2\pi\, i(   (\phi(x,y) -t)s    - x\cdot \xi  )\right)} $ is an eigenvector of the differential operator 
\begin{equation}\label{operator} 
\frac{1}{ 2\pi \, i(  s\frac{ \partial \phi(x,y)}{\partial x_l} -\xi_l)  }\frac{\partial }{\partial x_l}.
\end{equation}

Integrating by parts $M$ times,
we conclude that $$|G(s,\xi, \zeta)| \le c_M \sup_{(x,y)} \left|  s \frac{ \partial \phi(x,y)}{\partial x_l} - \xi_l     \right|^{-M}.$$

We interpret this estimate in three separate cases: $\left\{|\xi| \leq \frac{c}{2}|s|\right\}$, $\left\{2sC \leq |\xi|\right\}$, and $\left\{2sC >|\xi| > \frac{c}{2}|s|\right\}$.  
\\

If $|\xi| \leq \frac{c}{2}|s|$, then
 $$\left|  s \frac{ \partial \phi(x,y)}{\partial x_l} - \xi_l     \right|    \geq    \frac{1}{d} \left|   s\nabla_x\phi(x,y) -\xi  \right| \geq 
 \frac{1}{d} \left(  \left|   s\nabla_x\phi(x,y) \right|  -   \left|   \xi  \right| \right) \geq \frac{c}{2d}|s|.$$\\
 
Similarly, if  $2sC \leq |\xi|$, then  $$\left|  s \frac{ \partial \phi(x,y)}{\partial x_l} - \xi_l     \right| \geq   \frac{|\xi|}{2d}.$$

In either of the first two cases, if $|\xi|\sim 2^k$, $$|G(s, \xi, \zeta) | \le C_M \inf\{ |s|^{-M}, 2^{-kM}\},  $$
which in turn is $ \le C_M \inf\{ |s|^{-M}, 2^{-jM}, 2^{-kM}\} $ if $|\zeta|\sim 2^j \leq C|s|$.

In the case that $|\zeta|\sim 2^j >C|s|$, then the second inequality in \eqref{away} holds (otherwise $\frac{|\zeta|}{|s|}\in [c,C]$), and we can repeat the argument above to conclude that in this case, $$|G(s, \xi, \zeta) | \le C_M \inf\{ |s|^{-M}, 2^{-jM}, 2^{-kM}\}.  $$

Now, consider the case when $\{2sC  >|\xi| >\frac{c}{2}|s|\}$.  In this regime, we claim that the second equality in \eqref{critical} cannot hold.  Indeed, if the second equality in \eqref{critical} held, then $\{sC \geq |\zeta| \geq c|s|\}$. But then note that  if $|j-k| >K$, $|\zeta| \sim 2^j$, and $|\xi|\sim 2^k$, then both $\{2sC>|\xi| > \frac{c}{2}|s|\}$ and $\{sC>|\zeta| > c |s|\}$ cannot both simultaneously hold. 

Thus, when we can repeat the argument above (for the second inequality in \eqref{away}), and we are only left with the case that both $\{2sC>|\xi| > \frac{c}{2}|s|\}$ and $\{2sC>|\zeta| > \frac{c}{2}|s|\}$ hold. But again,  if $|j-k| >K$, $|\zeta| \sim 2^j$, and $|\xi|\sim 2^k$,  then both $\{2sC>|\xi| > \frac{c}{2}|s|\}$ and $\{2sC>|\zeta| > \frac{c}{2}|s|\}$ cannot both simultaneously hold.  \\

In all cases, we see that if $|j-k| >K$, $|\zeta| \sim 2^j$, and $|\xi|\sim 2^k$, then $$|G(s, \xi, \zeta) | \le C_M \inf\{ |s|^{-M}, 2^{-jM}, 2^{-kM}\}.  $$

\vskip.125in 

\subsection{Proof of Lemma \ref{fenergy}} 
The left hand side of \eqref{fenergyeq} can be re-written (see e.g. Proposition 8.5 in \cite{W04}) as
\begin{align*}
 \int {|\widehat{g\lambda}(\xi)|}^2 {|\xi|}^{-\gamma} d\xi 
&=c_{\gamma,d} \iint |x-y|^{\gamma -d }  g(x) g(y) d\lambda(y) d\lambda(x) \\
&=c_{\gamma,d} \int g(x)  Sg(x) d\lambda(x) \\
&\le c_{\gamma,d} \cdot \|g\|_{ L^2(\lambda)}\,   \cdot \|Sg\|_{ L^2(\lambda)}     \\
\end{align*}
where $Sg(x) = \int |x-y|^{\gamma-d}g(y) d\lambda(y) .$
We apply Schur's test (see Lemma 7.5 in \cite{W04}) to verify that there exists a constant $c>0$ so that if $s_{\lambda}> d-\gamma$, then
$$ \|Sg\|_{ L^2(\lambda)}  \le c  \|g\|_{ L^2(\lambda)} .$$
To complete the proof of the lemma, we need only verify that
Schur's test applies.  Observe that  
$$\int |x-y|^{\gamma-d}d\lambda(y)  =\int |x-y|^{\gamma-d}d\lambda(x),$$
 and assuming without loss of generality that the diameter of the support of $\lambda$ is at most $1$,

\begin{align*}
\int |x-y|^{\gamma-d}d\lambda(y) 
&= \sum_{j\geq 0 }   \int_{\left\{2^{-(j+1)}  \le |x-y| \le 2^{-j}  \right\}}  |x-y|^{\gamma -d} d\lambda(y)\\
&\le c' \sum_{j\geq 0 }  2^{j(d-\gamma)}\, \lambda(B(x,2^{-j})).
 \end{align*}
We recall that $   \lambda(B(x,2^{-j} )) \le  C 2^{-j\, s_{\lambda}}$ and conclude that 
the quantity above is bounded above when $s_{\lambda}> d-\gamma$.  
\vskip.125in

\section{Proof of Theorem \ref{pinnedtheoremchain}} 

As before, let $\mu$ be a Frostman measure supported on $E$. Recall that this means that 
\begin{equation} \label{muball} \mu(B(x,r)) \leq Cr^{s_{\mu}} \end{equation} for $s_{\mu}<dim_{{\mathcal H}}(E)$ and we can make $s_{\mu}$ as close to $dim_{{\mathcal H}}(E)$ as we want.
\vskip.125in 

Let $C_{k,\phi}^x(E)$ as in \eqref{kchain}, and
define the measure $\nu^{(k)}_x$ on $C_{k,\phi}^x(E)$ by the relation 
\begin{equation} \label{chainmeasdef} \int_{\mathbb{R}^k}f(t \,) d\nu^{(k)}_x(t \,)
=\int f\left(   \phi(x,x^2), \phi(x^2,x^3) , \cdots, \phi(x^{k},x^{k+1})   \right) d\mu(x^2)\cdots d\mu(x^{k+1}). \end{equation} 

Applying Cauchy-Schwartz we see that (formally) 
\begin{equation} \label{chaincsformal} 1  ={\left( \int d\nu^{(k)}_x(t) \right)}^2 \leq   |C_{k,\phi}^x(E)| \cdot \int  \left( \nu^{(k)}_x  (t) \right)^2 dt.   \end{equation} 
\vskip.125in 

As in the proof of Theorem \ref{generaltech}, the strategy is to establish the bound
$$ \int \int \left(\nu_x^{(k)}(t) \right)^2d\mu(x) dt \leq C,$$ and conclude that 
$$ \int  \left(\nu_x^{(k)}(t) \right)^2 dt<\infty$$ for $\mu$-a.e. $x \in E$.  Plugging this into (\ref{chaincsformal}) shows that $|\Delta^x_{\phi}(E)|>0$ for $\mu$ a.e. $x \in E$. 

In order to prove the more precise estimate (\ref{chainexceptionalset}) we will show that if $\lambda$ is a compactly supported Borel measure such that 
\begin{equation} \label{lambdaball} \lambda(B(x,r)) \leq Cr^{s_{\lambda}} \end{equation} for some $s_{\lambda}>0$, then 
\begin{equation}\label{chainkey} \int \int \left(\nu_x^{(k)}(t) \right)^2 d\lambda(x) dt \leq C \end{equation}
if $$s_{\lambda}+dim_{{\mathcal H}}(E)>d+1.$$ 

The conclusion (\ref{chainexceptionalset}) is recovered by taking $\lambda$ to be a Frostman measure on a subset of $E$. 
\vskip.125in 

We prove \eqref{chainkey} following the same basic reductions found in the proof of Theorem \ref{generaltech}. 
Let $\rho: \mathbb{R}^k \rightarrow \mathbb{R}^{+}$ be a smooth cut-off function, supported in the ball of radius $2$ and identically equal to $1$ in the ball of radius $1$ centered at the origin, with $\int \rho=1$. Let $\rho_{\epsilon}(u)=\epsilon^{-k} \rho(u/\epsilon)$. 
Then by the definition of the measure $\nu_x^{(k)}$ given in (\ref{chainmeasdef}), 
\begin{align*}
\nu_x^{(k)}*\rho_{\epsilon}(t)
&=\epsilon^{-k} \int \rho \left(\frac{t\,  -u}{\epsilon} \right) d\nu_x^{(k)}(u\,)\\
&=\,  \epsilon^{-k} \int  \rho \left(\frac{t\, -  \left(  \phi(x,x^2), \dots, \phi(x^k, x^{k+1})  \right)}{\epsilon} \right) d\mu(x^2) \cdots d\mu(x^{k+1}).
\end{align*}
\vskip.125in 

Let $\beta$ be a smooth cut-off function. Since the integration in $t$ is compact owing to the compactness of $E$, we can throw $\beta$ in for free. Squaring and integrating with respect to $d\lambda(x)$, we obtain 
$$ \int \int {\left| \frac{1}{\epsilon^k}  \int  \rho \left(\frac{t\, -  \left(  \phi(x,x^2), \dots, \phi(x^k, x^{k+1})  \right)}{\epsilon} \right) d\mu(x^2) \cdots d\mu(x^{k+1}) \right|}^2 d\lambda(x) \beta(t)dt.$$
\vskip.125in 

As in Theorem \ref{generaltech}, it follows that in order to bound $$\int \int \left(\nu_x^{(k)}(t) \right)^2 d\lambda(x) dt $$ 
it suffices to let $S^{\epsilon}_{\phi, t}$ equal
$$ \left\{ (x=x^1=y^1,x^2, \cdots, x^{k+1}, y^2, \cdots, y^{k+1}): \left| \phi(x^i,x^{i+1}) - t_i \right| \le \epsilon, \left| \phi(y^i,y^{i+1}) - t_i \right| \le 
\epsilon, 1\le i\le k \right\}$$ and show that 
\begin{equation}\label{long} \frac{1}{\epsilon^{2k}}\int \lambda \times \mu \times \cdots \times \mu \left\{ S^{\epsilon}_{\phi, t} \right\}  \beta(t) dt \le C,
\end{equation} independently of $\epsilon$.

\vskip.125in 

Letting $T^{\epsilon}_{\phi}$, the operator as in \eqref{operatordefn}, bounding the expression above reduces to bounding the following from above by a constant independent of $t$ in a compact set and $\epsilon>0$

\begin{equation}\label{comp}
\int \int {\left|  T^{\epsilon}_{\phi, t_1} \left(  T^{\epsilon}_{\phi, t_2} \left(  \cdots
\left(  T^{\epsilon}_{\phi, t_k} \left( T^{\epsilon}_{\phi, t_{k+1}} (\mu) \cdot \mu \right) \cdot \mu \right)  \right) \cdot \mu\right)\right|}^2 d\lambda(x) \beta(t)dt. \end{equation}
We repeatedly apply the mapping property in \eqref{opbound} $k$ times with $\lambda = \mu$ to see that this expression is bounded above when 
$s_{\lambda}>d+1-s_{\mu}$ and $s_{\mu}>\frac{d+1}{2}$.

\end{document}